\documentclass[11pt]{amsart}
\usepackage{amssymb}
\usepackage[arrow,matrix]{xy}
\usepackage{amsmath,amssymb, bbm, amscd, amsthm,mathrsfs,hyperref}
\theoremstyle{plain}
\newtheorem{thm}{Theorem}[section]
\newtheorem{cor}[thm]{Corollary}
\newtheorem{lem}[thm]{Lemma}
\newtheorem{prop}[thm]{Proposition}
\theoremstyle{definition}
\newtheorem{defn}[thm]{Definition}
\newtheorem{eg}[thm]{Example}

\newtheorem{rem}[thm]{Remark}

\def\al{\alpha}
\def\bt{\beta}
\def\dt{\delta}

\def\tt{\theta}
\def\vph{\varphi}
\def\lmd{\lambda}
\def\vps{\varepsilon}

\def\bd{\Delta}
\def\blmd{\Lambda}

\def\ms{\mathscr}
\def\mc{\mathcal}

\def\ot{\otimes}
\def\op{\oplus}

\def\se{\leqslant}
\def\le{\geqslant}

\def\Hom{\operatorname {Hom}}

\def\RHom{\operatorname {RHom}}
\def\Ext{\operatorname {Ext}}

\def\End{\operatorname {End}}

\def\dim{\operatorname {dim}}
\def\id{\operatorname {id}}

\def\Gr{\operatorname {Gr}}

\def\hdet{\operatorname {hdet}}
\def\GrMod{\operatorname {GrMod}}
\def\ev{\operatorname {ev}}
\def\db{\operatorname {coev}}

\def\t{\text}

\def\it{\textit}

\def\kk{\mathbbm{k}}
\def\ul{\underline}
\def\ra{\rightarrow}

\def\ZZ{\mathbb{Z}}
\def\NN{\mathbb{N}}

\begin{document}
\title{\bf Calabi-Yau Nichols algebras of Hecke type}

\author{Xiaolan Yu}
\address {Xiaolan YU\newline Hangzhou Normal University, Hangzhou, Zhejiang 310036, China}
\email{xlyu1005@yahoo.com.cn}

\author{Yinhuo Zhang}
\address {Yinhuo ZHANG\newline Department WNI, University of Hasselt, Universitaire Campus, 3590 Diepeenbeek,Belgium } \email{yinhuo.zhang@uhasselt.be}

\date{}

\begin{abstract}
Let $R$ be a  Nichols algebra of Hecke type. In this paper, we show that if $R$ is Noetherian of finite global dimension, then $R$ has a rigid dualizing complex.  We then give a necessary and sufficient condition for $R$ to be a Calabi-Yau algebra.
\end{abstract}

\keywords{Nichols algebra; Calabi-Yau algebra}
\subjclass[2000]{16T05,16E40. }

\maketitle

\section*{Introduction}

In \cite{yz2} we studied the Calabi-Yau (CY for short) property of pointed Hopf algebras and Nichols algebras of finite Cartan type. An interesting phenomenon is that the CY property of a pointed Hopf algebra $U(D,\lambda)$ of finite Cartan type is  dependent only on the CY property of its associated graded Hopf algebra $R\# kG=U(D,0)$, where $R$ is the associated Nichols algebra and $G$ is the coradical group of $U(D, \lambda)$. In other words,  if the smash product $U(D,0)$ is CY, then any lifting $U(D,\lambda)$ is CY. This raises a natural question whether a lifting or a cocycle deformation of a graded pointed CY Hopf algebra is still a CY algebra. At this point, we are not able to answer this question. On the other hand, we see that the graded pointed Hopf algebra $U(D,0)$ is a smash product of a Nichols algebra with a group algebra. One would expect that the CY property of $U(D,0)$ depends strongly on the CY property of the Nichols algebra $R$. However, this is not the case for the pointed Hopf algebras of finite Cartan type. In fact, when $U(D,\lambda)$ is CY, then $R$ must not be CY. Conversely if $R$ is CY, then $U(D,\lambda)$ is not CY. The reason was later discovered in \cite{yz} that the CY property of a smash product Hopf algebra $R\# H$ is strongly dependent on the action of $H$ on the algebra $R$, even when $H$ is a semisimple Hopf algebra.

In \cite[Thm 2.8]{yz}, we found that if $H$ is a semsimple Hopf algebra and $R$ is a CY Hopf algebra in the Yetter-Drinfeld module category over $H$, then $R\# H$ is CY if and only if the homological integral of $R$ is trivial. Thus we can construct CY Hopf algebras from CY braided Hopf algebras in the Yetter-Drinfeld module category over some (semsimple) Hopf algebra. On the other hand, the study of CY braided Hopf algebras has its own interest becasue it gives us examples of CY Hopf algebras in braided monoidal categories as we tend to work more and more with Hopf algebras in braided monoidal categories. This motivates us to study the CY property of certain important classes of braided Hopf algebras, e.g., Nichols algebras in the Yetter-Drinfeld module categories over group algebras as they play important role in the classification of pointed Hopf algebras. We did so in \cite{yz2}, where we obtained the rigid dualizing complex of a Nichols algebra of finite Cartan type and gave a necessary and sufficient condition for such a Nichols algebra to be Calabi-Yau.

In this paper we study the CY property of Nichols algebras of Hecke type. Unlike a Nichols algebra of finite Cartan type, a Nichols algebra of Hecke type,in general, cannot be expressed by generators and relations. However, we know that it is a quadratic algebra. Moreover, when its label $q$ is 1 or not a root of unity, it is a Koszul algebra \cite{as}.
Our main goal is to give a necessary and sufficient condition for a Nichols algebra of Hecke type to be Calabi-Yau. This is achieved in Section 2. The paper is organized as follows.

In Section 1, we recall some basic definitions and properties for braided Hopf algebras in the Yetter-Drinfeld module category over a Hopf algebra. We know that a finite dimensional braided Hopf algebra is Frobenus (see \cite{fms}). As a preparation for Section 2, we give the concrete formula of the Nakayama automorphism of a finite dimensional braided Hopf algebra (see Proposition \ref{naka}).

We prove our main theorem in Section \ref{2}. Let $V$ be a rigid braided vector space of Hecke type. We introduce the notion of a homological matrix of $V$.  Let  $R=\mc{B}(V)$ be the Nichols algebra of $V$. If $R$ is of finite global dimension, then we show that the rigid dualizing complex of $R$ is determined by its homological matrix and the braiding of $V$ (Theorem \ref{rid}). Hence, we are able to give a necessary and sufficient condition for the algebra $R$ to be a Calabi-Yau algebra.

\subsection*{Notations}
We work over a fixed algebraically closed field $\kk$ of characteristic zero. All algebras and vector spaces are assumed to be over $\kk$.

Given an algebra $A$, let $A^{op}$ denote the opposite algebra of
$A$ and $A^e$ denote the enveloping algebra $A\ot A^{op}$ of $A$. The unfurnished tensor $\ot$ means
over $\kk$ in this paper. Mod$A$ stands for the category of left $A$-modules.
We use Mod$A^{op}$ to denote the category of right $A$-modules.

A Noetherian algebra in this paper means a \textit{left and right}
Noetherian algebra.

For a Hopf algebra, we use Sweedler's notation (sumless) for its
comultiplication and its coactions.

\section{Finite dimensional braided Hopf algebras}\label{1}
\subsection{Braided vector spaces} We begin by recalling a braided vector space.  Let  $V$ be a vector space and $c : V\ot V \ra V\ot V$ a
linear isomorphism. Then $(V, c)$ is called a \textit{braided vector
space}, if $c$ is a solution of the following braid equation
$$ (c\ot \id )(\id\ot c)(c\ot \id ) = (\id \ot c)(c \ot \id )(\id\ot c).$$

Assume that $V$ is a finite dimensional vector space with a basis $\{v_1,\cdots,v_N\}$. Let $$\ev:V^*\ot V\ra \kk$$ and $$\db:\kk\ra V\ot V^*$$ be the evaluation map and coevaluation map respectively. That is, $$\ev(f\ot v)=f(v)\;\;\text{and}\;\;\db(1)=\sum_{i=1}^Nv_i\ot v_i^*.$$

A finite dimensional braided vector space $(V,c)$ is called \it{rigid} if the map $c^b:V^*\ot V\ra V\ot V^*$ defined by
$$c^b:=(\ev\ot \id_{V\ot V^*})(\id_{V^*}\ot c\ot \id_{V^*})(\id_{V^*\ot V}\ot \db)$$
is an isomorphism.

The notion of a braided vector space is closely related to the notion of a Yetter-Drinfeld module. Let $H$ be a Hopf algebra. A (left) Yetter-Drinfeld module $V$ over $H$ is simultaneously a left $H$-module and a left $H$-comodule satisfying the compatibility condition
$$\dt(h\cdot v)=h_1v_{(-1)}S_H(h_3)\ot h_2\cdot v_{(0)},$$ for any $v\in V$ and $h\in H$.

We denote by $^H_H\mc{YD}$ the category of Yetter-Drinfeld modules
over $H$ with  morphisms given by $H$-linear and $H$-colinear maps.

The tensor product of two Yetter-Drinfeld modules $M$ and $N$ is
again a Yetter-Drinfeld module with the module and the comodule
structures given as follows:
$$h\cdot(m\ot n)= h_1 \cdot m\ot  h_2\cdot n$$ $$\dt(m\ot  n) =
m_{(-1)}n_{(-1)}\ot  m_{(0)}\ot n_{(0)}, $$for any $h\in H$, $m\in
M$ and $n\in N$.  This turns the category of Yetter-Drinfeld modules
$^H_H\mc{YD}$ into a braided tensor category.

For any two Yetter-Drinfeld modules $M$ and $N$, the braiding
$c_{M,N} : M\ot N\ra  N \ot M $ is given by $$ c_{M,N}(m\ot n) =
m_{(-1)}\cdot n \ot m_{(0)},$$ for any $ m \in M$ and $n\in N$.
A braided vector space $(V,c)$ can be realized as a Yetter-Drinfeld module over some Hopf algebra $H$
if and only if $c$ is rigid \cite{ta}. If this is the case, it can be realized in many different ways. The FRT-construction is well-known \cite{frt}. The following theorem can be viewed as Hopf version of FRT-construction.

\begin{thm}\label{frt}\cite[Thm. 3.2.9]{sc} Let $(V, c)$ be a rigid braided vector space.
\begin{enumerate}
\item There exists a coquasitriangular Hopf algebra $H(c)$ such that V is a left $H(c)$-comodule and the braiding $c$ equals
the braiding on $V$ induced by the coquasitriangular structure of $H(c)$.
\item For all Hopf algebras $H$ having $V$ as a Yetter-Drinfeld module such that the induced braiding equals $c$,  there is a
unique morphism of Hopf algebras   $\psi: H(c)\ra H$ such that
$$\dt_H = (\psi\ot
 \id_M)\dt_{H(c)} \text{ and } u \cdot v= \psi (u) \cdot v,$$
 for any $u\in H(c)$ and $v\in V$.
\end{enumerate}
\end{thm}

\begin{rem}
Let $(H,\triangledown,u,\Delta,\vps,r)$ be a coquasitriangular Hopf algebra. For two $H$-comodules $M$ and $N$,
$c_{M,N} : M\ot N\ra  N \ot M $ given by $$ c_{M,N}(m\ot n) =
r(n_{(-1)}\ot m_{(-1)})n_{(0)} \ot m_{(0)},$$ for any $ m \in M$ and $n\in N$ defines a braiding. This turns the category of $H$-comodules to be a braided monoidal category.

Every $H$-comodule $M$ becomes a Yetter-Drinfeld module over $H$ with the action
$$h\cdot m=r(m_{(-1)}\ot h)m_{(0)},$$ for any $h\in H$ and $m\in M$.

\end{rem}
Let us recall the construction of $H(c)$. Let $\{v_1, v_2,\cdots, v_N\}$ be a basis of $V$, whose braiding is given by
 $$c(v_i\ot v_j)=\sum_{1\se m,n\se N}c_{ij}^{mn}v_m\ot v_n.$$The algebra $H(c)$ is generated by $T^i_j$ and $U^i_j$, $(i,j=1,\cdots,N)$ subject to the relations
$$\sum_{1\se k,l\se N}c_{ij}^{kl}T^m_kT^n_l-\sum_{1\se k,l\se N}T^k_iT^l_jc^{mn}_{kl};$$
$$\sum_{1\se k\se N}U^j_kT^k_i-\dt_{ij};$$
$$\sum_{1\se k\se N}T^j_kU^k_i-\dt_{ij}.$$

The Hopf algebra $H(c)$ coacts on $V$ as follows,
$$\dt(v_i)=\sum_{1\se j\se N}T^j_i\ot v_j.$$

\subsection{Braided Hopf algebras}
A braided Hopf aglebra is  of course defined as a Hopf algebra in a braided category.
\begin{defn}Let $H$ be a Hopf algebra.
\begin{enumerate}
\item[(i)] An \textit{algebra} in $^H_H\mc{YD}$ is a $\kk$-algebra $(R,m,u)$ such that
$R$ is a Yetter-Drinfeld $H$-module, and both the multiplication $m:R\otimes R\ra R$ and the unit $u:\kk\ra R$ are morphisms in $^H_H\mc{YD}$.
\item[(ii)] A   \textit{coalgebra}  in $^H_H\mc{YD}$ is
a $\kk$-coalgebra $(C,\Delta,\varepsilon)$ such that $C$ is a Yetter-Drinfeld $H$-module, and both the comultiplication $\bd:R\ra R\otimes R $ and the counit $\vps:R\ra \kk$  are morphisms in
$^H_H\mc{YD}$.
\end{enumerate}
\end{defn}

Let $R$ and $S$ be two algebras in $^H_H\mc{YD}$. Then $R\ot S$ is a
Yetter-Drinfeld module in $^H_H\mc{YD}$, and  becomes an algebra in
the category $^H_H\mc{YD}$ with  multiplication $m_{R\ul{\ot}S}$
given by $$m_{R\ul{\ot}S}=(m_R\ot m_S)(\id\ot c \ot \id).$$
Denote this algebra by $R\ul{\ot} S$.

\begin{defn}
Let $H$ be a Hopf algebra. A  \textit{braided bialgebra} in
$^H_H\mc{YD}$ is a 5-tuple $(R,m, u,\bd, \vps)$, where
\begin{enumerate}
\item[(i)] $(R,m, u)$ is an algebra in $^H_H\mc{YD}$.
\item[(ii)] $(R, \bd, \vps)$ is a coalgebra in $^H_H\mc{YD}$.
\item[(iii)] $\bd:R\ra R\ul\ot R$ and $\vps: R\ra \kk$ are  morphisms of algebras in $^H_H\mc{YD}$.
\end{enumerate}
If, in addition, the identity is convolution invertible in $\End
(R)$, then $R$ is called a  \textit{braided Hopf algebra} in
$^H_H\mc{YD}$. The inverse of the identity is called the
\textit{antipode} of $R$.
\end{defn}

In order to distinguish comultiplications of braided Hopf algebras
from those of usual Hopf algebras, we use Sweedler's notation with
upper indices for braided Hopf algebras
\begin{equation*}\label{sweed}\bd(r)=r^1\ot r^2\end{equation*}for $r\in R$.

Let $H$ be a Hopf algebra and $R$ a braided Hopf algebra in the
category  $^H_H\mc{YD}$.
The algebra $R\#H$ is a usual Hopf algebra with
the following structure \cite{ra}:

The multiplication is given by \begin{equation}(r\#g)(s\#h):=r(g_1\cdot s)\#g_2h\end{equation}
with unit $u_R\ot u_H$.
The comultiplication is given by
\begin{equation}\label{equa braidedcomulti}
\bd(r\#h):=r^1\#(r^2)_{(-1)}h_1\ot
(r^2)_{(0)}\#h_2
\end{equation}
with counit   $\vps_R\ot \vps_H$.
 The antipode is as follows:
\begin{equation}\label{equa braidedanti}S_{R\#H}(r\#h)=(1\#S_H(r_{(-1)}h))(S_R(r_{(0)})\#1).
\end{equation}

The algebra $R\#H$ is called the \textit{Radford biproduct} or the
\textit{bosonization} of $R$ by $H$. The algebra $R$ is a subalgebra of $R\#H$
and $H$ is a Hopf subalgebra of $R\#H$.

The following lemma is useful in this paper.

\begin{lem}\cite[Lemma 2.6]{yz}\label{2.6}
Let $H$ be a Hopf algebra, and $R$  a  braided Hopf algebra in the
category ${^H_H\mathcal{YD}}$. Then
$$ {S}_{R\#H}^2(r)={S_H(r_{(-1)})}\cdot {{S}_R^2(r_{(0)})},$$
for any $r\in R$.
\end{lem}

Now we have the following consequence.

\begin{cor}\label{inverse}
Let $H$ be a Hopf algebra, and $R$ a  braided Hopf algebra in the
category ${^H_H\mathcal{YD}}$. Then
$$ {S}_{R\#H}^{-2}(r)= r_{(-1)} \cdot {  {S}_R^2(r_{(0)})},$$
for any $r\in R$.
\end{cor}

\subsection{$\bt$-Frobenius extensions}

Now we recall some known results about $\bt$-Frobenius extensions.

Let $A$ be an algebra. For a left $A$-module $M$ and an algebra automorphism $\phi:A\ra A$, we denote by
$_{\phi}M$ the left $\phi$-twisted $A$-module with the underlying space $M$ and the left action $$a\cdot m=\phi(a)m$$ for any $a\in A$ and $m\in M$.  Similarly, for a right $A$-module $N$, we have
$N_{\phi}$.

\begin{defn}
Let $B\subseteq A$ be an algebra extension and $\bt:B\ra B$ an algebra automorphism. The extension $B\subseteq A$ is called a (left) \textit{$\bt$-Frobenius extension} if \begin{enumerate}
\item $A$ is a finitely generated projective right $B$-module;
\item $A\cong {}_\bt\Hom_B(A_B, B_B)$ as $B$-$A$-bimodules, where  ${}_\bt\Hom_B(A_B, B_B)$ is a left $\bt$-twisted $B$-$A$-bimodule via $(b\cdot \psi\cdot a)(x)=\bt(b)\psi(ax)$, for all $b\in B$, $a, x\in A$, and $\psi\in \Hom_B(A_B, B_B)$.
\end{enumerate}
\end{defn}

Note that in the case where $\bt=id$ and $B =\kk$, the algebra $A$ is a Frobenius algebra. The isomorphism
$A_A\ra (A^*)_A$ of Part (2) in the definition is the classical Frobenius isomorphism.

\begin{lem}\cite[Thm. 1.1]{bf},\cite[Prop. 1.3]{fms}
Let $B\subseteq A$ be an algebra extension, $\bt: B\ra B$ an algebra automorphism, and $f : A\ra {}_\bt B$ a $B$-$B$-bimodule map. The following are equivalent:
\begin{enumerate}
\item $B\subseteq A$ is a $\bt$-Frobenius extension;
\item There exist $r_i, l_i\in A$, $i = 1,\cdots, n$, such that for any $a\in A$,
\begin{enumerate}
\item $a=\sum_{i=1}^nr_if(l_ia);$
\item $a=\sum_{i=1}^n (\bt^{-1}f)(ar_i)l_i.$
\end{enumerate}
\end{enumerate}
\end{lem}

A morphism $f:A\ra B$ satisfying part (2) is called a \textit{$\bt$-Frobenius homomorphism} and  $\{r_i, l_i\}$  is called a \textit{dual pair basis} for $B\subseteq A$. In the classical case of a Frobenius algebra $A$, one can define the Nakayama automorphism $\eta$ of $A$ by $f(xy)=f(y\eta(x))$, for all $x$, $y\in A$.

\subsection{Finite dimensional braided Hopf algebras}
One can define (left and right) integrals for any augmented algebras. Let $A$ be an augmented algebra with augmentation $\vps:A\ra \kk$. Then $\int^l_A=\{t\in A|at=\vps(a)t\}$ is the space of left integrals. The space of right integrals can be defined similarly.

Let $R$ be a finite dimensional braided Hopf algebra in the category $^H_H\mc{YD}$. Let $A=R\# H$ and  $\overline{A}=A/AH^+$. Then by \cite[Lemma 5.5]{fms}, the inclusion $R\subseteq A$ induces an isomorphism of coalgebras $R\cong \overline{A}$. The following lemma follows from Theorem 3.1 and Corollary 3.3 in \cite{fms}.

\begin{lem}\label{frob}
Let $H$ be a Hopf algebra, $R$ a finite dimensional braided Hopf algebra in the category $^H_H\mc{YD}$ and $A=R\#H$. Then $\int^r_{R^*}=\int^r_{\overline{A}^*}=\kk\lmd$ for some $\lmd\neq 0$. We may choose a $\blmd\in R$ such that $\lmd\leftharpoonup\blmd=\vps$, and there exists an algebra homomorphism $\chi:H\ra \kk$ such that $$\lmd(h\cdot r)=\chi(h)\lmd(r),$$for all $r\in R$ and $h\in H$.

Define $$f:A\ra H \;\;\t{by}\;\; f(a)=\lmd(\overline{a_1})a_2$$
$$\bt:H\ra H \;\;\t{by}\;\; \bt(h)=\chi(h_1)h_2$$
for all $a\in A$, $h\in H$.  Then:
\begin{enumerate}
\item $H\subseteq A$ is a $\bt$-Frobenius extension with Frobenius map $f$.
\item $R$ is a Frobenius algebra with Frobenius homomorphism $f_R$, where $f_R$ is the restriction of $f$ to $R$.
\item There exist $x_i, y_i\in R$, $1\se i\se m$, such that in $A{}_\bt\ot_H A$, $$\sum x_i\ot y_i=\sum S^{-1}(\blmd_2)\ot \blmd_1.$$
The set $\{S^{-1}(\blmd_2), \blmd_1\}$ forms a dual pair basis for $H\subseteq A$, and $\{x_i,y_i\}$ gives a dual pair basis for $R$.
\end{enumerate}

\end{lem}

It is easy to see that $\blmd$ in Lemma \ref{frob} belongs to $\int_R^r$. By Corollary 5.8 in \cite{fms}, we know that $\int_R^r$ is one-dimensional. Hence, we can define the \textit{right modular function} $\al\in R^*$  for $R$ as  $r\blmd=\al(r)\blmd$, for all $r\in R$.

We have seen that a finite dimensional braided Hopf algebra is Frobenius. The following proposition delivers the Nakayama automorphism.
\begin{prop}\label{naka}
Let $H$ be a Hopf algebra and $R$ a finite dimensional braided Hopf algebra in the category $^H_H\mc{YD}$.  Choose $0\neq\blmd\in \int_R^r$. Assume that there exists an algebra homomorphism $\chi_{_H}: H\ra \kk$ such that $h\cdot\blmd=\chi_{_H}(h)\blmd$. Then the Nakayama automorphism $\eta:R\ra R$ is defined as follows:
$$\begin{array}{ccl}\eta(r)&=&\chi_{_H}((r^2)_{(-1)})\al(r^1)S_A^{-2}((r^2)_{(0)})\\&=&\chi_{_H}((r^2)_{(-1)1})\al(r^1)(r^2)_{(-1)2}\cdot S_R^{-2}((r^2)_{(0)}).\end{array}$$
\end{prop}

\proof Let $A=R\#H$ and keep the same notations as in Lemma \ref{frob}. Since $\{S^{-1}(\blmd_2),\blmd_1\}$ forms a dual pair basis and there exist $x_i, y_i\in R$, such that $\sum x_i\ot y_i=\sum S^{-1}(\blmd_2)\ot \blmd_1$ in $A{}_\bt\ot_H A$, we have the following:
$$\begin{array}{rcl}\eta(r)&=&S^{-1}(\blmd_2)f(\blmd_1\eta(r))\\
&=&x_if(y_i\eta(r))\\
&=&x_if(ry_i),
\end{array}$$since the restriction of $f$ to $R$ is the Frobenius homomorphism of $R$. Thus
$$\begin{array}{ccl}S_A^2(\eta(r))&=&f(ry_i)S^2_A(x_i)\\
&=&f(r\blmd_1)S_A(\blmd_2)\\
&=&\lmd(\overline{r_1\blmd_1})r_2\blmd_2S_A(\blmd_3)\\
&=&\lmd(\overline{r_1\blmd})r_2\\
&=&\lmd(\overline{r^1\#(r^2)_{(-1)}\blmd})(r^2)_{(0)}\\
&=&\lmd(\overline{r^1((r^2)_{(-1)1}\cdot\blmd)\#(r^2)_{(-1)2}})(r^2)_{(0)}\\
&=&\chi_{_H}((r^2)_{(-1)1})\vps((r^2)_{(-1)2})\lmd(r^1\blmd)(r^2)_{(0)}\\
&=&\chi_{_H}((r^2)_{(-1)})\al(r^1)\lmd(\blmd)(r^2)_{(0)}\\
&=&\chi_{_H}((r^2)_{(-1)})\al(r^1)(r^2)_{(0)}.
\end{array}$$
Therefore, $\eta(r)=\chi_{_H}((r^2)_{(-1)})\al(r^1)S_A^{-2}(r^2)_{(0)}$. By Corollary \ref{inverse}, we have $S_A^{-2}(r)=r_{(-1)}\cdot S_R^{-2}(r_{(0)}),$ for any $r\in R$.

Hence, $$\begin{array}{ccl}\eta(r)&=&\chi_{_H}((r^2)_{(-1)})\al(r^1)S_A^{-2}((r^2)_{(0)})\\&=&\chi_{_H}((r^2)_{(-1)1})\al(r^1)(r^2)_{(-1)2}\cdot S_R^{-2}((r^2)_{(0)}).\end{array}$$ \qed

\section{The Calabi-Yau property}\label{2}
\subsection{Graded algebras} An $\NN$-graded algebra $R=\op_{i\le0}R_i$ is called \textit{connected} if $R_0=\kk$. In the sequel, a graded algebra always means a \textit{connected graded} algebra generated in degree 1. We denote by $\GrMod(R)$ the category of graded left
$R$-modules with graded homomorphisms of degree 0. For graded modules $M$ and $N$, let $M(n)$ be the shifted module of $M$ with $M(n)_i=M_{n+i}$, and
$$\underline{\Hom}_R(M,N)=\bigoplus_{i\in \ZZ}\Hom_{\GrMod(R)}(M,N(i)).$$
Let  $\underline{\Ext}_R^i(-,-)$ be the derived functor of $\underline{\Hom}_R(-,-)$.  If $M$ is finitely generated, then
$\Hom_R(M,N)=\underline{\Hom}_R(M,N)$. If $R$ is in addition Noetherian, then
$\Ext^i_R(M,N)=\underline{\Ext}_R^i(M,N)$ for all $i\le 0$.

Now, we recall some basic concepts related to quadratic algebras and Koszul algebras. Let $V$ be a finite dimensional vector space, and $T(V)$ the tensor algebra with the
usual grading. A graded algebra $R=T(V)/(I)$ is called a \textit{quadratic} algebra if $I$ is a
subspace of $V\ot V$. The quadratic dual of $R$ is defined as the algebra $R^! = T(V^*)/(I^\perp)$,
where $I^\perp$ is the orthogonal subspace of $I$ in $V\ot V$.

For convenience, we identify $(V^*)^{\ot n}$ with $(V^{\ot n})^*$ by means of the following:
$$(f_n\ot f_{n-1}\ot\cdots f_1)(x_1\ot x_2\ot\cdots x_n)=f_1(x_1)f_2(x_2)\cdots f_n(x_n).$$

\begin{defn} A quadratic algebra R is called \textit{Koszul} if the trivial R-module $_R\kk$
admits a graded projective resolution
$$\cdots\ra P_n\ra P_{n-1}\ra\cdots\ra P_1\ra P_0\ra _R\kk\ra 0 $$
such that $P_n$ is generated in degree $n$ for all $n\le 0$.
\end{defn}

The Koszul complex of a quadratic algebra is important in the study of quadratic algebras. Let $\{v_i\}_{i=1,2,\cdots, N}$ be a basis of $V$ and  $\{v_i^*\}_{i=1,2,\cdots,N}$ the dual basis of $V^*$.  The complex $\mc{K}_l(R)$ is the complex
$$\cdots\ra R\ot R^{!*}_m\ra  R\ot R^{!*}_{m-1}\ra \cdots R\ot R_2^{!*}\ra R\ot V \ra R.$$
The differential $d_m:R\ot R^{!*}_m\ra  R\ot R^{!*}_{m-1}$ is defined by
$d_m(a\ot \al)=\sum_{i=1}^N av_i\ot v_i^*\cdot\al $, where $v_i^*\cdot\al $ is defined by $(v_i^*\cdot\al )(x)=\al(xv_i^*)$ for any $x\in A_{m-1}^!$.
\begin{lem}\cite[Thm. 2.6.1]{bgs}
Let $R=T(V)/(I)$ be a quadratic algebra. Then the following are equivalent:
\begin{enumerate}\item $R$ is a Koszul algebra;
\item The complex $\mc{K}_l(R)$ is a resolution of $\kk$.
\end{enumerate}
\end{lem}

We refer to \cite{bgs} for more details about Koszul algebras and the Koszul duality.

\begin{defn}
A graded algebra $R$ is called \textit{AS-Gorenstein} (Artin-Schelter Gorenstein) with parameters $(d,l)$ for some integers $d$ and $l$, if
\begin{enumerate}
\item both injective dimensions of $_RR$ and of $R_R$ are equal to $d$;
\item $\underline{\Ext}^i_R(\kk,R)\cong \underline{\Ext}^i_{R^{op}}(\kk,R)\cong \begin{cases}0 ,& i\neq d ;
\\ \kk(l),& i=d .\end{cases}$
\end{enumerate}
If, furthermore, $R$ is of finite global dimension, then $R$ is called \textit{AS-regular}.
\end{defn}

\subsection{Homological determinants}Let $H$ be a Hopf algebra. We call a graded algebra $R$  a \textit{graded $H$-module algebra}, if $R$ is a left $H$-module algebra and each component $R_i$ is also a left $H$-submodule.
 The concept of a homological determinant  was defined by J{\o}rgensen and Zhang \cite{joz} for graded automorphisms of an
AS-Gorenstein algebra and by Kirkman, Kuzmanovich and
Zhang \cite{kkz} for Hopf actions on an AS-Gorenstein algebra.

\begin{defn}\label{hdetl}(cf.\cite{kkz})
Let $H$ be a Hopf algebra, and  $R$ an $H$-module AS-Gorenstein algebra with injective dimension $d$.  There is a natural  $H$-action on $\underline{\Ext}_R^d(\kk,R)$ induced by the
$H$-action on $R$. Let $\bf{e}$ be a non-zero element in
$\underline{\Ext}_R^d(\kk,R)$. Then there exists an algebra homomorphism $\eta:
H\ra\kk$ satisfying $h \cdot {\bf e} = \eta(h)\bf{e}$ for all $h\in
H$. The composite map $\eta S_H:H\ra\kk$  is called the \textit{
homological determinant} of the $H$-action on $R$, denoted
by $\hdet $ (or more precisely $\hdet_R$).
\end{defn}

Next, we give the homological determinant of a Koszul algebra. It is probably well known. Since we do not find any reference, we give a proof here for completeness.  Let $R$ be a quadratic $H$-module algebra and $R^!$ its quadratic dual.  Then $R_1=V$ is a left $H$-module. Hence, $V^*$ is a right  $H$-module with the $H$-action defined by \begin{equation}\label{triangle}(f\triangleleft h)(v)=f(h\cdot v),\end{equation}
for any $f\in V^*$ and $v\in V$. We can make $V^*$ into a left $H$-module by using the antipode:
\begin{equation}\label{ac0}h\cdot f=(f\triangleleft S_H^{-1})(h),
\end{equation} for any $f\in V^*$ and $h\in H$. The action ``$\triangleleft$'' can be extend to $(V^*)^{\ot n}$ by $$(f_n\ot f_{n-1}\ot \cdots \ot f_1)\triangleleft h=f_n\triangleleft {h_n}\ot f_{n-1}\triangleleft h_{n-1}\ot\cdots \ot f_1\triangleleft h_1.$$ It is easy to check that
this defines an $H$-module structure on $R^!$.

\begin{lem}\label{H-Koszul}
Let $H$ be a Hopf algebra and $R$ a left graded $H$-module algebra. If $R$ is Koszul AS-regular, then the $R$-projective resolution of $\mc{K}_l(R)\ra\kk\ra 0$ of $\kk$ is an $R\#H$-module complex.
\end{lem}
\proof For each $m\in \NN$, $R^!_m$ is a right $H$-module, so $R^{!*}_m$ is a left $H$-module with the action
$$(h\cdot \al)(x) =\al(x\triangleleft h)$$
for any $h\in H$, $\al \in A^{!*}_m $ and $x\in R^{!}_m$. For any $x\#h\in R\#H$ and $r\ot \al\in R\ot R^{!*}_m$, let
$$(x\#h)\cdot (r\ot \al)=x (h_1\cdot r)\ot h_2\cdot \al.$$
Now each $R\ot R^{!*}_m$ is a left $R\#H$-module. To show that the differentials in $\mc{K}_l(R)\ra\kk\ra 0$ are morphisms of $R\#H$-morphisms, we only need to check that they are $H$-morphisms. Suppose that $h\cdot v_i=\sum_{j=1}^Nc_{ji}^hv_j$ for $h\in H$ and $1\se i\se n$. Then $v_i^*\triangleleft h=\sum _{j=1}^N c^h_{ij}v^*_j$.

We then have
$$\begin{array}{ccl}h\cdot (d_m(r\ot \al))&=&h\cdot (\sum_{i=1}^Nrv_i\ot  v^*_i\cdot\al)\\
&=&\sum_{i=1}^Nh_1\cdot(rv_i)\ot h_2\cdot( v^*_i\cdot\al)\\
&=&\sum_{i=1}^N(h_1\cdot r )(h_2\cdot v_i)\ot h_3\cdot(v^*_i\cdot\al)\\
&=&\sum_{i=1}^N\sum_{j=1}^N(h_1\cdot r)c_{ji}^{h_2}v_j\ot h_3\cdot(v^*_i\cdot\al)\\
&=&\sum_{i=1}^N\sum_{j=1}^N(h_1\cdot r)v_j\ot c_{ji}^{h_2}h_3\cdot(v^*_i\cdot\al).
\end{array}$$
Now for any $x\in R^!_m$, we have $$\begin{array}{ccl}\sum_{i=1}^Nc_{ji}^{h_1}h_2\cdot(v^*_i\cdot\al)(x)&=&\sum_{i=1}^Nc_{ji}^{h_1}(v^*_i\cdot\al)(x\triangleleft h_2)\\
&=&\sum_{i=1}^Nc_{ji}^{h_1}\al((x\triangleleft h_2) v^*_i)\\
&=&\sum_{i=1}^N\al( (x\triangleleft h_2)c_{ji}^{h_1}v^*_i)\\
&=&\al( (x\triangleleft h_2)(v^*_j\triangleleft h_1))\\
&=&\al( (xv^*_j)\triangleleft h)\\
&=&h\cdot \al(xv^*_j)\\
&=&( v^*_j\cdot(h\cdot \al))(x).
\end{array}$$
So, $$\begin{array}{ccl}\sum_{i=1}^N\sum_{j=1}^N(h_1\cdot r)v_j\ot c_{ji}^{h_2}h_3\cdot(v^*_i\cdot\al)&=&\sum_{j=1}^N(h_1\cdot r)v_j\ot v^*_j\cdot(h_2\cdot \al) \\
&=&d_m (h\cdot(r\ot \al)).
\end{array}$$Therefore, the differentials are $H$-morphisms. We have completed the proof.\qed

\begin{prop}
Let $R$ be a Noetherian left graded $H$-module algebra. Assume that $R$ is Koszul AS-regular with global dimension $d$. Then $R^!_d=\kk\blmd$ for some $\blmd\neq 0$. The homological determinant of $R$ is given by $\hdet(h)\blmd=\blmd \triangleleft h$ for any $h\in H$.
\end{prop}

\proof  Since $R$ is Koszul, AS-regular and of global dimension $d$, by \cite[Prop. 5.10]{sm} $R^!$ is Frobenius with $\dim R^!_d=1$. Choose a nonzero element $\blmd$ in $ R^!_d$. So $R^!_d=\kk\blmd$.

 Applying $\underline{\Hom}_R(-,R)$ to the resolution $\mc{K}_l(R)$ of $\kk$, we obtain that $\underline{\Ext}_R^d(\kk,R)$ is the cohomology at the final position of the following complex of left $H$- and right $R$-modules (not necessary $H$-$R$-bimodules):
\begin{equation}\label{com1}0\ra \underline{\Hom}_R(R,R)\ra \underline{\Hom}_R(R\ot R_2^{!*},R)\ra\cdots \ra \underline{\Hom}_R(R\ot R_d^{!*},R)\ra0.\end{equation}

Note that each $R^!_m$ is finite dimensional. So we have $\underline{\Hom}_R(R\ot R_m^{!*},R)\cong R_m^!\ot R$. Therefore, the complex (\ref{com1}) is isomorphic to the following complex of left $H$- and right $R$-modules:
$$0\ra R\ra V^*\ot R\ra R_2^!\ot R\ra \cdots\ra R_d^!\ot R\ra 0,$$
where for each $1\se m\se d$, the $H$-action on $R_m^!\ot R$ is given by $h\cdot (f\ot r)=h_1\cdot f\ot h_2\cdot r$ for any $h\in H$, $f\in R_m^!$ and $r\in R$. The nonzero element $\blmd$ in $R^!_d$ sits in the $d$-th cohomology. Otherwise the $d$-th cohomology is zero. The 1-dimensional component $R^!_d$ is an $H$-module. Now for any $h\in H$,  on one hand, we have $$h \cdot \blmd =\hdet(S^{-1}_H(h))\blmd.$$ On the other hand, $$h \cdot \blmd =\blmd \triangleleft S_H^{-1}(h).$$ Therefore, for any $h\in H$, $$\hdet(h)\blmd=\blmd \triangleleft h. $$ \qed
\subsection{Calabi-Yau algebras}
 We follow Ginzburg's definition of a Calabi-Yau algebra \cite{g2}.

\begin{defn}A graded algebra $A$ is called a  \textit{graded Calabi-Yau algebra of dimension
$d$} if
\begin{enumerate}
\item[(i)] $A$ is \it{homologically smooth}, that is, $A$, as a graded $A^e$-module, has
a bounded resolution of finitely generated projective
$A$-$A$-bimodules;
\item[(ii)] There are $A$-$A$-bimodule
isomorphisms$$\Ext_{A^e}^i(A,A^e)=\begin{cases}0& i\neq d
\\A(l)&i=d.\end{cases}$$ as graded $A^e$-modules for some integer $l$.
\end{enumerate}
\end{defn}
In the sequel, Calabi-Yau will be abbreviated to CY for short.

CY algebras form a class of algebras possessing a rigid dualizing complex. The
non-commutative version of a dualizing complex was first introduced
by Yekutieli.

\begin{defn}\cite{y} (or \cite[Defn. 6.1]{vdb})\label{defn dc}
Assume that $A$ is a Noetherian algebra. Then an object
$\ms{R}$ of $D^b(A^e)$ is called a \textit{dualizing complex}  if it satisfies the
following conditions:
\begin{enumerate}
\item[(i)] $\ms{R}$ is of  finite injective dimension over $A$ and $A^{op}$.
\item[(ii)] The cohomology of $\ms{R}$ is given by bimodules which are
finitely generated on both sides.
\item[(iii)] The natural morphisms $A\ra  {\RHom}_A(\ms{R},\ms{R})$ and $A\ra
\RHom_{A^{op}}(\ms{R},\ms{R})$ are isomorphisms in $D(A^e)$.
\end{enumerate}
\end{defn}

Roughly speaking, a dualizing complex is a complex $\ms{R}\in
D^b(A^e)$ such that the functor
$$\RHom_A(-,\ms{R}):D^b_{fg}(A)\ra D^b_{fg}(A^{op})$$ is a
duality, with adjoint $\RHom_{A^{op}}(-,\ms{R})$ (cf. \cite[Prop.
3.4 and Prop. 3.5]{y}). Here $D^b_{fg}(A)$ is the  full triangulated
subcategory of $D(A)$ consisting of bounded complexes with finitely
generated cohomology modules.

Dualizing complexes are not unique  up to isomorphism. To overcome
this weakness,  Van den Bergh introduced the concept of a rigid
dualizing complex cf. \cite[Defn. 8.1]{vdb}.

\begin{defn}\label{defn rdc} Let $A$ be a Noetherian (graded) algebra. A dualizing  complex $\ms{R}$ over $A$ is called \textit{rigid} (in the graded case) if $$\RHom_{A^e}(A,{_A\ms{R}\ot \ms{R}_A})\cong
 \ms{R}$$ in $D(A^e)$ ($D(\Gr A^e)$).
\end{defn}

The following lemma follows from \cite[Prop. 8.2 and Prop. 8.4]{vdb}.

\begin{lem}\label{cor cyrid}
Let $A$ be a Noetherian graded algebra which is  homologically smooth. Then
$A$ is a graded CY algebra of dimension $d$ if and only if $A$ has a rigid
dualizing complex $A[d](-l)$ for some $d, l\in \ZZ$.
\end{lem}

\subsection{Nichols algebras} In this paper, we will focus on the Calabi-Yau property of Nichols algebras of Hecke type. Nichols algebras form an important class of braided Hopf algebras generated by primitive
elements. They play important roles in the classification of pointed Hopf algebras. They appeared  first in the paper \cite{ni} of Nichols.

\begin{defn}
Let $V$ be a Yetter-Drinfeld module over a Hopf algebra $H$. A
graded braided Hopf algebra $R=\op_{i\le 0}R_i$ in the category
$^H_H\mc{YD}$ is called a  \textit{Nichols algebra} of $V$ if the
following conditions hold:
\begin{enumerate}
\item $R_0\cong \kk$ and $R_1\cong V$.
\item $R_1=P(R)$, the set of primitive elements in $R$.
\item $R$ is generated  as an algebra by $R_1$.
\end{enumerate}
We denote by $\mc{B}(V)$ the Nichols algebra $R$.
\end{defn}

\begin{defn} We say that a braided vector space $(V, c)$ is of \textit{Hecke-type} with label $q\in \kk$,
$q\neq 0$, if
$$(c - q)(c + 1) = 0.$$
\end{defn}

In the following, we assume that the label $q$ is either 1 or not a root of 1.

\begin{lem}\cite[Prop. 3.4]{as}\label{nichols hecke}
Let $(V,c)$ be a braided vector space of Hecke-type with label $q$, which is either 1 or not a root of 1. Then $\mc{B}(V)$ is a quadratic algebra.
Moreover, $\mc{B}(V)$ is a Koszul algebra and its quadratic dual is the Nichols algebra $\mc{B}(V^*)$ corresponding to the braided vector space $(V^*, -q^{-1}c^t)$.
\end{lem}

\subsection{Main results}

We need some preparations to prove the main theorem.

\begin{lem}\label{braidiongVstar}
Let $V$ be a finite dimensional Yetter-Drinfeld module over a Hopf algebra $H$, $c$ the braiding of $V$ induced by the Yetter-Drinfeld module structure. Define a Yetter-Drinfeld module on $V^*$ as follows
\begin{equation}\label{ac}(h\centerdot f)(v)=-q^{-1}(h\cdot f)(v)=-q^{-1}f(S^{-1}_H(h)\cdot v),
\end{equation}
for all $h\in H, f\in V^*, v\in V$, where $h\cdot f$ is given by $\mathrm{(\ref{ac0})}$ . For any $f\in V^*$, $\dt(f)=f_{(-1)}\ot f_{(0)},$ is determined by the equation \begin{equation}\label{coa} f_{(-1)}f_{(0)}(v)=S_{H}(v_{(-1)})f(v_{(0)}), \end{equation} for any $ v\in V$.
The braiding of $V^*$ induced by the Yetter-Drinfeld module structure is just $-q^{-1}c^t$.
\end{lem}
\proof It is straightforward to check that the action and the coaction defined in the lemma give a Yetter-Drinfeld module structure on $V^*$. Next we show that the braiding induced by the Yetter-Drinfeld module structure is just $-q^{-1}c^t$.  Let $\{v_1,v_2,\cdots,v_N\}$ be a basis of $V$. Assume that for each $1\se i\se N$, $$\dt(v_i)=h^{i1}\ot v_1+h^{i2}\ot v_2+\cdots+h^{iN}\ot v_N,$$ with $h^{i1},h^{i2},\cdots, h^{iN}\in H$.  Further, we assume that for each $1\se i\se N$, $1\se j\se N$, $$( h^{ij} \cdot v_1, h^{ij} \cdot v_2,\cdots, h^{ij}\cdot v_N)=(v_1,v_2,\cdots,v_N)H^{ij},$$ where $H^{ij}$ is an $N\times N$ matrix over $\kk$. Then the braiding induced by the Yetter-Drinfeld module structure is given by
\begin{equation}\label{braidiongV}c(v_i\ot v_j)=\sum_{m=1}^N\sum_{n=1}^N H^{in}_{mj}v_m\ot v_n.\end{equation}
The transposed braiding $c^t$ of $V^*$ is given by \begin{equation}\label{braidiongV*}c^t(v_i^*\ot v_j^*)=\sum_{k=1}^N\sum_{l=1}^N H^{mi}_{jn}v_n^*\ot v_m^*.\end{equation} By equation (\ref{coa}), we obtain that
$$\dt(v_i^*)=S_H(h^{1i})\ot v_1^*+S_H(h^{2i})\ot v_2^*+\cdots+S_H(h^{Ni})\ot v_N^*,$$for each $1\se i\se N$. Therefore, the braiding $\overline{c}$ of $V^*$ induced by the Yetter-Drinfeld module structure is given by $$\begin{array}{ccl}\overline{c}(v_i^*\ot v_j^*)&=&v^*_{i(-1)}\centerdot v^*_j \ot v^*_{i(0)}\\&=&\sum_{m=1}^NS_H(h^{mi})\centerdot v_j^*\ot v^*_m\\
&=&-q^{-1}\sum_{l=1}^N\sum_{k=1}^N H^{mi}_{jn} v_n^*\ot v^*_m.
\end{array}$$
It follows from Equation (\ref{braidiongV*}) that  $\overline{c}$ coincides with the braiding $-q^{-1}c^t$. \qed

Now we can see that the Nichols algebra $\mc{B}(V)$ is AS-regular and Koszul.

\begin{prop}\label{AS} Let $V$ be  a rigid braided vector space of hecke type with label $q$, which is either 1 or not a root of 1. Let $R=\mc{B}(V)$ be the Nichols algebra of $V$. If $R$ is a Noetherian algebra of global dimension $d$, then
\begin{enumerate}
\item The algebra $R$ is a Koszul algebra;
\item The algebra $R^!$ is Frobenius;
\item The algebra  $R$ is AS-regular.
\end{enumerate}
\end{prop}

\proof The fact that the algebra $R$ is Koszul follows directly from Lemma \ref{nichols hecke}.

Now let $c$ be the braiding of $V$. Since $V$ is rigid, the vector space $V$ is a Yetter-Drinfeld module over the algebra $H(c)$ by Theorem \ref{frt}. Furthermore, since the global dimension of $R$ is $d$,  the quadratic dual $R^!=\op_{i=0}^{d}R^!_i$ is finite dimensional. Hence by Lemma \ref{nichols hecke} and Lemma \ref{braidiongVstar}, the algebra $R^!=\mc{B}(V^*)$ is a finite dimensional braided Hopf algebra in the category $^{H(c)}_{H(c)}\mc{YD}$, where the Yetter-Drinfeld module structure of $V^*$ is given by equations (\ref{ac}) and (\ref{coa}).  It follows from Lemma \ref{frob} that $R^!$ is Frobenius.

Now $R^!$ is a Frobenius algebra.  By \cite[Prop. 5.10]{sm}, the algebra $R$ is AS-Gorenstein. Moreover, the global dimension of $R$ is finite. Thus $R$ is AS-regular. \qed

Now the Nichols algebra $R$ is a  Koszul AS-regular algebra with global dimension $d$. Hence, we have $\dim R^!_d=1$. That is, $R^!_d=\kk\blmd$ for some nonzero element $\blmd$ in $ R^!_d$. The linear map $V^*\ra
V^*$ sending any $v^*\in V^*$ to $-q^{-1}v^*$ extends to a graded automorphism $\psi$ of $R^!$. The automorphism $\psi$ restricts to a linear map of the 1-dimensional component $R^!_d$. Therefore, there exists $Q\in \kk$, such that $\psi(\blmd)=Q\blmd$.  We call $Q$ the \textit{quantum label} of $R$.

Recall that the algebra $H(c)$ is generated by $T^i_j$ and $U^i_j$ $(1\se i,j \se N)$, where $N$ is the dimension of the braided vector space $V$. Let $d_{ij}=\hdet(T^i_j)$.  We call the matrix $D=\{d_{ij}\}$ the \textit{homological matrix} of $R$.

Let $R=T(V)/(I)$ be a quadratic algebra, and $R^!$ the quadratic dual of $R$. Assume that  $\vph$ is an automorphism on $R$ and $\psi$ is  an automorphism on $R^!$. We say that $\vph$ and $\psi$ are dual to each other if the two restrictions $\vph|_V$ and $\psi|_{V^*}$  are dual to each other as automorphisms of vector spaces.

\begin{thm}\label{rid} Let $V$ be a rigid braided vector space of Hecke-type with label $q$, which is either 1 or not a root of 1. Suppose that the braiding is given as follows,
$$c(v_i\ot v_j)=\sum_{1\se m,n\se N}c^{mn}_{ij}v_m\ot v_n,$$
where $N$ is the dimension of $V$ and $v_1,v_2,\cdots, v_N$ is a basis of $V$. Let $R=\mc{B}(V)$ be the Nichols algebra of $V$. If $R$ is a Noetherian algebra of global dimension $d$, then
\begin{enumerate}
\item The rigid dualizing complex of  $R$ is isomorphic to $ _{\phi\epsilon^{n+1}}R[d](-d)$, where $\epsilon$ is the automorphism of $R$ given by the multiplication by $(-1)^m$ on $R_m$, and $\phi$ is the automorphism defined by $$\phi(v_i)=-q^{-1}Q\sum^N_{l=1}\sum^N_{j=1}\sum^N_{k=1}d_{lk}c_{ji}^{jk}v_l,$$ for each $1\se i\se N$, where $Q$ and $D=\{d_{ij}\}$ are the quantum label and the homological matrix of $R$ respectively.
\item $R$ is a CY algebra if and only if for each $1\se i\se N$,
 $$-q^{-1}Q\sum^N_{j=1}\sum^N_{l=1}d_{lk}c_{ji}^{jk}=\begin{cases}(-1)^{d+1} &  l=i ;
\\ 0&  l\neq i.\end{cases}$$
\end{enumerate}
\end{thm}

\proof (1) Let $H=H(c)$. We first clarify the $H$-action on $V$ and $V^*$. Since the braiding is given as follows,
$$c(v_i\ot v_j)=\sum_{1\se m,n\se N}c^{mn}_{ij}v_m\ot v_n.$$ for each $1\se i,j,n\se N$.  We obtain that
$$T^n_i\cdot v_j=\sum_{m=1}^Nc^{mn}_{ij}v_m,$$ for each $1\se i,j,n\se N$. By Equation (\ref{ac}), we have
 \begin{equation}\label{3}S_H(T^k_j)\centerdot v_j^*=\sum_{l=1}^N-qc_{jl}^{jk}v_l^*,
 \end{equation}
for each $1\se k,j\se N$.

From the proof of Proposition \ref{AS}, we see that the algebra $\mc{B}(V^*)$ is a finite dimensional braided Hopf algebra in the category $^H_H\mc{YD}$. The Yetter-Drinfeld module structure of $V^*$ is given by Equations (\ref{ac}) and (\ref{coa}). Hence, the algebra $R^!$ is Frobenius. Moreover, the algebra $R$ is Koszul, AS-regular, and of global dimension $d$. So $R^!_d=\kk\blmd$ for some $\blmd\neq 0$. It is easy to see that $\blmd$ belongs to $\int_{R^!}^r$ as well. Following from Proposition \ref{H-Koszul}, we obtain $$\hdet(h)\blmd=\blmd\triangleleft h,$$
for any $h\in H$. This yields the following: \begin{equation}\label{chi}h\centerdot \blmd=Qh\cdot \blmd=Q\blmd\triangleleft S_H^{-1}(h)=Q\hdet(S_H^{-1}(h))\blmd,\end{equation}
where the first equation follows from (\ref{ac}). It follows that the assumptions of Proposition \ref{naka} hold.

Now let $A=R\#H$ and $A^!=R^!\#H$. From the proof of Lemma \ref{braidiongVstar}, we see that for each $1\se i\se N$,
$$\dt(v_i^*)=S_H(T^i_1)\ot v_1^*+S_H(T^i_2)\ot v_2^*+\cdots+S_H(T^i_N)\ot v_N^*.$$

By Corollary \ref{inverse}, we have that $$\begin{array}{ccl}S^{-2}_{A^!}(v_i^*)&=&(v_i^*)_{(-1)}\centerdot (v_i^*)_{(0)}\\&=&\sum_{j=1}^NS_H(T^i_j)\centerdot v_j^*.\end{array}$$

Since $R^!=\mc{B}(V^*)$, the elements of $V^*$ are primitive. That is, for each $v^*\in V^*$, we have $\Delta_{R^!}(v^*)=v^*\ot 1+1\ot v^*$. Now following from Proposition \ref{naka}, we obtain that the Nakayama automorphism of $R^!$ is given by
$$\begin{array}{ccl}\eta(v_i^*)&\overset{\tiny{(\ref{chi})}}=&Q\sum_{j=1}^N\sum_{k=1}^N\hdet(S_H^{-1}((v_i^*)_{(-1)})S^{-2}_{A^!}((v_i^*)_{(0)})\\&=&Q\sum_{j=1}^N\sum_{k=1}^N\hdet(T^i_k)S_H(T^k_j)\centerdot v_j^*\\
&\overset{\tiny{(\ref{3})}}=&-q^{-1}Q\sum_{j=1}^N\sum_{k=1}^N\sum_{l=1}^Nd_{ik} c^{jk}_{jl}v^*_l.
\end{array}$$
The dual automorphism of this automorphism is given by
$$\phi(v_i)=-q^{-1}Q\sum^N_{l=1}\sum^N_{j=1}\sum^N_{k=1}d_{lk}c_{ji}^{jk}v_l.$$

Therefore, by \cite[Thm. 9.2]{vdb}, we obtain that the rigid dualizing complex of $R$ is isomorphic to $_{\phi\epsilon^{d+1}}R[d](-d)$. Note that our expression about the rigid dualizing complex is a bit different from the one in \cite[Thm. 9.2]{vdb}. This is caused by the identification of $(V^*)^{\ot n}=(V^{\ot n})^*$ in a different way.

(2) Since $R$ is a Noetherian Koszul AS-regular algebra, it is homologically smooth. Thus Part (2) follows from Part (1) and Corollary \ref{cor cyrid}.\qed

We end this section with two examples.

\begin{eg} Let $(V, c)$ be a braided vector space with $c$ given by  $$c(v_i\ot v_j)=q_{ij}v_j\ot v_i$$
for each $1\se i,j\se n$, where $\{v_1,v_2\cdots,v_n\}$ is a basis of $V$. We assume that $q_{ii}=1$ for $1\se i\se n$ and  $q_{ij}'s$ satisfy  $q_{ij}q_{ji}=1$ for each $1\se i,j\se n$. Then $(V,c)$ is a braided vector space of Hecke type with label 1. The Nichols algebra $R=\mc{B}(V)$ is isomorphic to the following algebra $$ \kk\langle v_1, \cdots ,  v_n\mid v_iv_j=q_{ij}v_jv_i, 1 \se i, j \se \tt, i\neq
j\rangle.$$ It is easy to check that $R$ is a Koszul algebra of global dimension $n$. The quadratic dual $R^!$ is isomorphic to $$ \kk\langle v^*_1, \cdots ,  v^*_n\mid v^*_iv^*_j=-q_{ij}v^*_jv^*_i, 1 \se i, j \se n, i\neq
j\rangle.$$
It is not difficult to see  $R^!_n=\kk\blmd$, where $\blmd=v_1^*v^*_2\cdots v^*_n$. So the quantum label $Q$ is $(-1)^n$ and the homological matrix is $$\left(\begin{array}{cccccc}q_{12}q_{13}\cdots q_{1n}&0&0&\cdots&0&0\\0&q_{21}q_{23}\cdots q_{2n}&0&\cdots&0&0\\\cdots&\cdots&\cdots&\cdots&\cdots&\cdots\\0&0&0&\cdots &0&q_{n1}q_{n2}\cdots q_{n,n-1}\end{array}\right).$$  By Theorem \ref{rid}, we obtain that the rigid dualizing
complex of $R$ is isomorphic to $ _\phi R(-n)[n]$, where $\phi$ is
the algebra automorphism defined by
$$\phi(v_i)=q_{i1}\cdots
q_{i(i-1)}q_{i(i+1)}\cdots q_{in}v_i,$$ for all $1\se i\se
n$. The algebra $R$ is CY if and only if $q_{i1}\cdots
q_{i(i-1)}q_{i(i+1)}\cdots q_{in}=1$ for all $1\se i\se n$.

\end{eg}

\begin{eg} Let $(V, c)$ be a braided vector space with basis $\{v_1,v_2,v_3,v_4\}$. The matrix of the braiding $c$ with respect to the basis $\{v_1\ot v_1,\cdots,v_1\ot v_4,v_2\ot v_1,\cdots,v_2\ot v_4,v_3\ot v_1,\cdots,v_3\ot v_4,v_4\ot v_1,\cdots,v_4\ot v_4\}$ is given by
$$\left(\begin{array}{cccccccccccccccc}
1&0&0&0&0&0&0&0&0&0&0&0&0&0&0&0\\
0&0&0&0&0&0&0&0&0&0&0&1&0&0&0&0\\
0&0&0&0&0&0&0&1&0&0&0&0&0&0&0&0\\
0&0&0&0&0&0&0&0&0&0&0&0&1&0&0&0\\

0&0&0&0&0&0&0&0&0&0&0&0&0&0&1&0\\
0&0&0&0&0&1&0&0&0&0&0&0&0&0&0&0\\
0&0&0&0&0&0&0&0&0&1&0&0&0&0&0&0\\
0&0&1&0&0&0&0&0&0&0&0&0&0&0&0&0\\

0&0&0&0&0&0&0&0&0&0&0&0&0&1&0&0\\
0&0&0&0&0&0&1&0&0&0&0&0&0&0&0&0\\
0&0&0&0&0&0&0&0&0&0&1&0&0&0&0&0\\
0&1&0&0&0&0&0&0&0&0&0&0&0&0&0&0\\

0&0&0&1&0&0&0&0&0&0&0&0&0&0&0&0\\
0&0&0&0&0&0&0&0&1&0&0&0&0&0&0&0\\
0&0&0&0&1&0&0&0&0&0&0&0&0&0&0&0\\
0&0&0&0&0&0&0&0&0&0&0&0&0&0&0&1\end{array}\right).$$
The braided space $(V,c)$ is of Hecke type with label 1. The Nichols algebra $R=\mc{B}(V)$ is isomorphic to the algebra generated by $v_1,v_2,v_3,v_4$ subject to the following relations:
$$v_1v_2=v_3v_4,\;\;\;v_1v_3=v_2v_4,\;\;\;v_4v_2=v_3v_1$$
$$v_4v_3=v_2v_1,\;\;\;v_1v_4=v_4v_1,\;\;\;v_2v_3=v_3v_2.$$
The algebra $R$ is a Koszul algebra of global dimension 4.  Let $R^!$ be the quadratic dual. We have $R^!_4=\kk\blmd$, where $\blmd=v_1^*v^*_2v^*_3 v^*_4$. So the quantum label $Q$ is 1 and the homological matrix is $$\left(\begin{array}{cccc}1&0&0&0\\0&1&0&0\\0&0&1&0\\0&0&0&1\end{array}\right).$$  By Theorem \ref{rid},  the rigid dualizing
complex of $R$ is isomorphic to $R[4](-4)$. Therefore, the algebra $R$ is a CY algebra of dimension 4.
\end{eg}

\vspace{5mm}

\bibliography{}

\end{document}